\documentclass[preprint,12pt]{elsarticle}

\usepackage{graphicx}%
\usepackage{multirow}%
\usepackage{amsmath,amssymb,amsfonts}%
\usepackage{amsthm}%
\usepackage{mathrsfs}%
\usepackage[title]{appendix}%
\usepackage{xcolor}%
\usepackage{textcomp}%
\usepackage{manyfoot}%
\usepackage{booktabs}%
\usepackage{algorithm}%
\usepackage{algorithmicx}%
\usepackage{algpseudocode}%
\usepackage{listings}%

\theoremstyle{plain}
\newtheorem{theorem}{Theorem}[section]

\newtheorem{corollary}[theorem]{Corollary}
\newtheorem{proposition}[theorem]{Proposition}

\theoremstyle{definition}
\newtheorem{definition}[theorem]{Definition}
\newtheorem{example}[theorem]{Example}
\newtheorem{Lemma}[theorem]{Lemma}
\theoremstyle{remark}
\newtheorem{remark}{Remark}

\newcommand \s{^{*}}
\newcommand \+{^{\dag}}
\newcommand \p{^{\perp}}
\newcommand \q{^{\times}}

\begin{document}
\begin{frontmatter}
    

\title{On unbounded complementable operators}


\author[inst1]{Sachin Manjunath Naik}

\affiliation[inst1]{organization={Department of Mathematics},
            addressline={Manipal Institute of Technology,\\ Manipal Academy of Higher Education}, 
            city={Manipal},
            postcode={576104}, 
            state={Karnataka},
            country={India}}

\author[inst2]{P. Sam Johnson}
\affiliation[inst2]{organization={Department of Mathematical and Computational Sciences},
            addressline={\\National Institute of Technology Karnataka, Surathkal}, 
            city={Mangaluru},
            postcode={575025}, 
            state={Karnataka},
            country={India}}

\begin{abstract}
The concept of complementability is extended from bounded operators to densely defined operators on Hilbert spaces. By introducing appropriate projections and decomposition techniques, a framework is developed for analyzing complementability in this broader context. The results provide new insights into the structure of unbounded operators, contributing to the ongoing development of operator theory.
\end{abstract}




\end{frontmatter}




\section{Introduction} 
 ~~ \newline
 The study of complementability in operator theory has traditionally focused on bounded linear operators on Hilbert spaces. However, many operators of significant interest, particularly those arising in quantum mechanics, partial differential equations, and spectral theory, are unbounded. Extending the concept of complementability to unbounded operators introduces new analytical challenges and offers potential applications in both pure and applied mathematics.

Unbounded operators $A: \mathcal{D}(A) \subseteq \mathcal{H} \to \mathcal{K}$ between Hilbert spaces $\mathcal{H}$ and $\mathcal{K}$ are often difficult to handle, primarily because their domain $\mathcal{D}(A)$ is generally a proper subspace of $\mathcal{H}$ and may fail to be densely defined or closed. Consequently, the notion of complementability must be redefined to accommodate these concepts. In particular, establishing the existence of projections that suitably decompose the domain and range of $A$ necessitates a more delicate treatment than in the bounded case.

A natural starting point for this investigation is the extension of the Schur complement to unbounded operators, inspired by the foundational work of many mathematicians. The Schur complement has undergone extensive development since its inception, leading to various generalizations and extensions. Ando \cite{Ando} introduced an alternative definition of the Schur complement relative to a subspace of $\mathbb{C}^n$, utilizing orthogonal projections to analyze the interaction between the matrix and associated operators. In particular, Ando introduced the concept of complementable operators with respect to a subspace in finite-dimensional Hilbert spaces, elucidating how these projections interact with the matrix and related operators.

Mitra and Puri \cite{Mitra} extended Ando's framework to $m \times n$ matrices over $\mathbb{C}$ by considering subspaces $M \subseteq \mathbb{C}^m$ and $N \subseteq \mathbb{C}^n$. Their work demonstrated that this generalization holds over arbitrary fields. Independently, Carlson \cite{CARLSON} further expanded the Schur complement to arbitrary fields, significantly broadening its scope.

More recently, Antezana et al. \cite{Antezana} extended these concepts to bounded linear operators on Hilbert spaces. They introduced notions of complementability and weak complementability for a bounded linear operator $T: \mathcal{H} \to \mathcal{K}$ with respect to two closed subspaces $M \subset \mathcal{H}$ and $N \subset \mathcal{K}$. This extension facilitated the definition of a generalized Schur complement for operators on Hilbert spaces. Further developments on complementable operators can be found in \cite{Albert1, Arias, Ben-Israel2003, Butler, Carlson1, CARLSON, Corach1, Corach2, Cottle, Mitra, SachinLAA}.

The implications of complementability theory are extensive. In spectral theory, it provides insights into the spectral decomposition of unbounded operators. In the study of partial differential operators, it plays a role in understanding boundary conditions and extensions. Furthermore, in control theory, unbounded operators naturally arise in infinite-dimensional dynamical systems, where complementability considerations are essential.

This work aims to establish a rigorous framework for complementability in the unbounded setting, investigate its fundamental properties, and explore its applications. By generalizing established notions from bounded operators to unbounded ones, we seek to bridge existing gaps in the literature and pave the way for further advancements in operator theory and its applications.

\section{Preliminaries}
~~ \newline
Let $\mathcal H$ and $\mathcal K $ be real or complex Hilbert spaces.  We denote  $\mathcal L(\mathcal H, \mathcal K),$  the space of all linear operators from $\mathcal H$ to $\mathcal K.$ We denote  $\mathcal B(\mathcal H, \mathcal K),$  the space of all linear bounded operators from $\mathcal H$ to $\mathcal K$ and we abbreviate $\mathcal B(\mathcal H)=\mathcal B(\mathcal H, \mathcal H)$. For $T \in \mathcal B(\mathcal H, \mathcal K),$ we denote by $T^*$, $\mathcal N(T)$ and $\mathcal R(T),$ respectively, the adjoint, the null-space and the range-space of $T$. For a subspace $M \subseteq \mathcal{H},$  we denote the set $\{x \in  M : \|x\| \leq 1 \}$ as $\mathcal{B}_M$. An operator $T\in \mathcal B(\mathcal H)$ is called a projection if $T^2=T.$ A projection, $T$ is called orthogonal if $T=T^*.$  An operator $T\in \mathcal B(\mathcal H)$ is called self-adjoint if $T\s=T.$ A a self-adjoint operator is called positive (denoted $T\geq 0$) if  $\langle Tx, x \rangle \geq 0$ for all $x\in \mathcal H$.

\begin{definition}(\cite{Antezana})
	Let $P_r \in \mathcal{B}(\mathcal{H})$ and $P_\ell \in \mathcal{B}(\mathcal{K})$ be projections. An operator $T \in \mathcal {B}(\mathcal{H}, \mathcal{K})$ is called $(P_r, P_\ell)$-complementable if there exist operators $M_r \in \mathcal{B}(\mathcal{H})$ and $M_\ell \in \mathcal{B}(\mathcal{K})$ such that 
	\begin{enumerate}
		\item $(I-P_r)M_r=M_r$ and $(I-P_\ell)TM_r=(I-P_\ell) T$,
		\item $M_\ell (I-P_\ell)=M_\ell$ and $M_\ell T(I-P_r)=T(I-P_r) $.
	\end{enumerate}
\end{definition}

This indicates that the definition is independent of the specific projections considered and depends solely on the ranges of these projections. The following theorem shows that this property depends only on the ranges of the projections.
 \begin{theorem}(\cite{Antezana})
\label{propdef}
	Let $P_r \in \mathcal{B}(\mathcal{H})$ and $P_\ell \in \mathcal{B}(\mathcal{K})$ be two projections whose ranges are  $M$ and $N$ respectively. Given $T \in \mathcal {B}(\mathcal{H}, \mathcal{K}),$ the following statements are equivalent:  
	\begin{enumerate}
		\item $T$ is $(P_r, P_\ell)$-complementable;
		\item $\mathcal{R}(C)\subseteq \mathcal{R}(D)$ and $\mathcal{R}(B\s) \subseteq \mathcal{R}(D\s)$;
		\item There exist two projections $E \in \mathcal{B}(\mathcal{H})$ and $F \in  \mathcal{B}(\mathcal{K})$ such that $\mathcal{R}(E ^*)=M$, $\mathcal{R}(F )=N$,  $\mathcal{R}(T E )\subseteq N$ and $\mathcal{R}( (FT)\s )\subseteq M. $ 
	\end{enumerate}
 \end{theorem}

 \begin{definition}[\cite{Antezana}]
	An operator $T \in \mathcal {B}(\mathcal{H}, \mathcal{K})$ is called $(M, N)$-complementable if it is $(P_r, P_\ell)$-complementable, for some projections $P_r $ and $ P_\ell$ with $\mathcal{R}(P_r)=M$ and $\mathcal{R}(P_\ell)=N$.
\end{definition}
\noindent If $\mathcal{H}= \mathcal{K}$, we use the notation $M$-complementable for $(M, M)$-complementable.

The following theorem, known as Douglas' Lemma or Douglas' Factorization Theorem, provides conditions under which a bounded linear operator can be factored through another operator.

\begin{theorem}[\cite{Douglas}]\label{dgls}
		Let $A, B \in \mathcal{B}(\mathcal H)$.  Then the following statements are equivalent:
		\begin{enumerate}
			\item $\mathcal{R}(A) \subseteq \mathcal{R}(B)$;
			\item $AA\s \leq \lambda BB\s$, for some $\lambda >0$;
			\item There exists a bounded operator $C \in \mathcal{B}(\mathcal H)$ such that $A=BC$.
		\end{enumerate}
		Moreover, if these equivalent conditions hold, then there is a unique operator $C \in \mathcal{B}(\mathcal H)$ such that
		\begin{enumerate}
			\item[(i)] $\|C\| = \inf \{ \lambda >0: AA\s \leq \lambda BB\s\} $;
			\item[(ii)] $\mathcal{N}(A) = \mathcal{N}(B)$;
			\item[(iii)] $\mathcal{R}(C) \subseteq \mathcal{N}(B)\p$.
		\end{enumerate}
	\end{theorem}
\noindent	The unique solution $C$ above is referred to as the Douglas reduced solution (or simply, the reduced solution) \cite{redsoln}.

Antezana et al.\cite{Antezana} extended the notion of the Schur complement from matrices to bounded operators on Hilbert spaces. Employing the Douglas factorization lemma, they defined the Schur complement of an operator relative to closed subspaces of Hilbert spaces. This formulation generalizes the classical Schur complement and provides a framework for analyzing more intricate operator structures in Hilbert spaces.
	\begin{definition}[\cite{IJPAM}]
		Let $T$  be $(M,N)$-complementable. Let $Z$ and $Y$ be the reduced solutions of the equations $C=DX$ and $B\s= D\s X$ respectively. Then, the Schur complement (the bilateral shorted operator) of $T$ with respect to closed subspaces 
		$M$ and $N$ is given by   $$ \begin{pmatrix}
			A-BZ & 0 \\
			0 & 0
		\end{pmatrix}=\begin{pmatrix}
     A-YC & 0 \\
     0 & 0
 \end{pmatrix}$$ and is denoted by $T_{/(M,N)}$.\\ If $\mathcal{H}=\mathcal K$ and $M=N$, then we denote $T_{/(M,N)}$ by $T_{/M}$. 
\end{definition}

\begin{theorem}[\cite{Antezana}]\label{C}
		Let $T$ be $(M,N)$-complementable. Then the following statements hold good :
		\begin{enumerate}
			\item $T\s$ is $(N,M)$-complementable;
			\item $\mathcal{R}(T_{/(M,N)})=\mathcal{R}(T) \cap N$;
			\item $\mathcal{N}(T_{/(M,N)})= M\p + \mathcal{N}(T)$;
			\item 	If $\mathcal{R}(D)$ is closed, then $T_{/(M,N)}=\begin{pmatrix}
				A-BD\+ C & 0\\
				0 & 0
			\end{pmatrix}.$
		\end{enumerate}
	\end{theorem}

    \begin{theorem}\label{eqgm}(\cite{IJPAM})
	Let $T \in \mathcal{B}(\mathcal{H},\mathcal{K})$ with $\mathcal{R}(D)$ a closed subspace of $\mathcal{K}$. Then $T$ is $(M,N)$-complementable if and only if for each $x \in M$, there exists a unique $z\in M$ such that $\{T(x,0)+T(M\p)\} \cap N = \{(z,0)\}$. Moreover, $T_{/(M,N)}(x,0)=(z,0)$.
	\end{theorem}

An operator $T : \mathcal{D}(T) \subseteq \mathcal {H} \to \mathcal {K}$ is called \textit{unbounded} if it is not bounded, i.e., for each $ \alpha >0$ there exists  $x \in \mathcal{D}(T) $  such that $ \|T\|>\alpha \|x\|.$
As a result, specifying a domain is a crucial part of defining an unbounded operator.

\begin{definition} [\cite{Kreyszig}] Let  $T ,U  \in \mathcal L(\mathcal {H}, \mathcal {K}),V \in \mathcal L( \mathcal {K},\mathcal {G})$
	and $\alpha \in \mathbb C \backslash \{0\}$. Then
	\begin{enumerate}
		\item $T + U \in \mathcal L(\mathcal {H}, \mathcal {K})$ with domain $\mathcal{D}(T + U) := {\mathcal{D}(T) \cap \mathcal{D}(U)}$ and
		$	(T + U)x = Tx + Ux$
		
		for all $x \in \mathcal{D}(T + U).$
		\item $VT \in \mathcal L(\mathcal {H}, \mathcal {G})$ with domain $\mathcal{D}(VT) = \{x \in \mathcal{D}(T) : Tx \in \mathcal{D}(V)\}$ and
		$(VT)x = V(Tx)$
		for all $x \in \mathcal{D}(VT).$
		\item $\alpha T \in \mathcal L(\mathcal {H}, \mathcal {K})$ with domain $\mathcal{D}(\alpha T) = \mathcal{D}(T)$ and
		$(\alpha T)x = \alpha Tx$
		for all $x \in \mathcal{D}(T).$
	\end{enumerate}

\end{definition}

\begin{definition} [\cite{Kreyszig}] 
Let $T$ be an operator such that an operator $T : \mathcal{D}(T) \subseteq \mathcal {H} \to \mathcal {K}$. \begin{enumerate}
    \item $T$ is called densely defined if $\mathcal{D}(T)$ is dense in $\mathcal {H}$. 
    \item $T$ is called closed if its graph $\mathcal{G}(T)=\{(x,Tx): x \in \mathcal{D}(T) \}$ is closed in $\mathcal {H} \times \mathcal {K}$.
\end{enumerate}
\end{definition}

\begin{theorem}\label{Douglasunbdd}(\cite{Douglas})
   Let $T:\mathcal{D}(T) \subseteq \mathcal{H}_1 \to \mathcal{K}$ and $S:\mathcal{D}(S) \subseteq \mathcal{H}_2 \to \mathcal{K}$ be operators satisfying $\mathcal{R}(T) \subseteq \mathcal{R}(S)$. Then there exists an operator $U:\mathcal{D}(T) \to \mathcal{D}(S)$ such that $T = SU$.
    \end{theorem}

\begin{definition}[\cite{Kreyszig}] Let $T$ and $U $ be operators from $\mathcal {H}$ to $\mathcal {K}$ such that $\mathcal{D}(T) \subseteq \mathcal{D}(U)$ and
$Tx = Ux$
for all $x \in \mathcal{D}(T).$ Then $U$ is called an extension of $T$ (or $T$ is called a restriction of $U$) and is
denoted by $U\supseteq T$ (or by $T \subseteq U$). Note that $T = U$ if and only if $T \subseteq U$  and  $U \subseteq T.$
\end{definition}

\begin{definition} [\cite{Kreyszig}]\label{adjoint_dfn} Let $T : \mathcal{D}(T) \subseteq \mathcal {H} \to \mathcal {K}$ be a densely defined operator.
Then, there exists a unique operator $T^*$ such that
	$$\langle T x, y\rangle = \langle x, T^* y\rangle \quad \text{ for all }
x \in \mathcal{D}(T ) \text{ and } y \in \mathcal{D}(T^* ).$$ This operator is called the adjoint of $T$. Here,
	$$   \mathcal{D}(T^*)=\Big\{y\in   \mathcal {K}: x\mapsto \langle Tx,y\rangle \mbox{ for all } x\in   \mathcal{D}(T) \mbox{ is continuous}\Big\}.$$
In other words,
$$
\mathcal{D}(T^* ) = \Big\{y \in \mathcal {K} :\text{  there exists a unique } z\in \mathcal {H},\langle T x, y\rangle = \langle x, z\rangle 
\text{ for all }x \in \mathcal{D}(T)\Big\}$$
and in this case,
$T^*y =z.$
\end{definition} 

\begin{proposition}(\cite{Kreyszig})

\begin{enumerate}

	\item  If $T \in \mathcal L (\mathcal {H}, \mathcal {K})$ is densely defined, then  	$(\alpha T)^* = \overline{\alpha} T^*,$   ~for any scalar $\alpha.$

	\item  If $T, U$ are densely defined such that $\mathcal{D}(TU)$ is dense, then $(TU)^* \supseteq U^* T^*.$
	If $T$ is everywhere defined, then $(TU)^* = U^* T^*$.
	\item If $T$ is densely defined such that $T \subseteq U,$ then $U^* \subseteq T^* .$
\end{enumerate}
\end{proposition}

\section{Operator decomposable subspaces}
In the study of (possibly unbounded) linear operators between Hilbert spaces, the structure of the operator’s domain plays a fundamental role in understanding various analytical and geometric properties, such as decomposability, closability, and the behavior of restrictions and extensions. Given a closed subspace \( M \subseteq \mathcal{H} \), it is natural to ask whether the domain \( \mathcal{D}(T) \) of an operator \( T: \mathcal{D}(T) \subseteq \mathcal{H} \to \mathcal{K} \) admits a decomposition compatible with the orthogonal decomposition \( \mathcal{H} = M \oplus M^{\perp} \). This leads to the notion of \(T\)-decomposability of a subspace \( M \), which requires that \( \mathcal{D}(T) \) splits as the direct sum of its intersections with \( M \) and \( M^{\perp} \). It is worth noting that this property holds automatically when \( T \) is bounded, since in that case \( \mathcal{D}(T) = \mathcal{H} \). However, for unbounded operators, \(T\)-decomposability becomes a nontrivial condition and provides useful insight into how the domain of \( T \) interacts with the subspace structure of \( \mathcal{H} \). The following definition formalizes this concept and sets the stage for further analysis.
\begin{definition}
Let \( T: \mathcal{D}(T) \subseteq \mathcal{H} \to \mathcal{K} \) be an operator and \( M \) be a closed subspace of \(\mathcal{H}\). We say that \( M \) is \( T \)-decomposable if
$$
\mathcal{D}(T) = (M \cap \mathcal{D}(T)) \oplus (M\p \cap \mathcal{D}(T)).
$$
\end{definition}

 \noindent The subspaces \( M \cap \mathcal{D}(T) \) and \( M\p \cap \mathcal{D}(T) \) are subspaces of \( \mathcal{D}(T) \). Hence, the inclusion 
$$
\mathcal{D}(T) \supseteq (M \cap \mathcal{D}(T)) \oplus (M\p \cap \mathcal{D}(T))
$$
always holds. Therefore, to show a subspace \( M \) to be \( T \)-decomposable, it is enough to show that 
$$
\mathcal{D}(T) \subseteq (M \cap \mathcal{D}(T)) \oplus (M\p \cap \mathcal{D}(T)).
$$

\noindent We use notations \(M_1= M \cap \mathcal{D}(T) \)  and \(M_2= M\p \cap \mathcal{D}(T) \).

\begin{remark}
   \begin{enumerate}
    \item If $M$ is any closed subspace of $\mathcal{D}(T)$, then $M$ is $T$-decomposable. Indeed, let $x \in \mathcal{D}(T)$. Then $x=y+z$ where $y \in M$ and $z \in M\p$.
     Since $M \subseteq \mathcal{D}(T)$, we get $$y \in M \cap \mathcal{D}(T) \text{~and~} z = x-y \in M\p \cap \mathcal{D}(T).$$ 
     Hence $$\mathcal{D}(T) = (M \cap \mathcal{D}(T)) \oplus (M\p \cap \mathcal{D}(T)).$$
 \item   If $\mathcal{N}(T)$ is closed, $\mathcal{N}(T)$ is $T$-decomposable.
\end{enumerate} 
\end{remark}
\noindent \begin{example}
Consider the operator $ T: \mathcal{D}(T) \subseteq \ell_2 \to \ell_2 $ defined by

\[
T(x_1, x_2, \ldots) = (x_1 - x_2, 0, 3(x_3 - x_4), 0, \ldots),
\]
where the domain \(\mathcal{D}(T)\) is given by

\[
\mathcal{D}(T) = \left\{ (x_n) : \sum_{k=1}^{\infty}(2k-1)^2 |x_{2k-1} - x_{2k}|^2  < \infty \right\}.
\]

\noindent Let \( M = \{ (x_n) : x_{2k} = 0 \text{ for all } k \} \). We can see that \( M \) is a closed subspace and $$ M\p = \{ (x_n) : x_{2k-1} = 0 \text{ for all } k \}.$$
Now consider the sequence $\Big( \frac{1}{n} \Big )$. We have
\[
\sum_{k=1}^{\infty} \left|\frac{1}{2k-1} - \frac{1}{2k}\right|^2 (2k-1)^2 = \sum_{k=1}^{\infty} \frac{1}{4k^2} < \infty,
\]
so \( \Big(\frac{1}{n}\Big) \in \mathcal{D}(T) \).
However, we can decompose \( \Big(\frac{1}{n}\Big) \) as

\[
\left(\frac{1}{n}\right) = \left(1, 0, \frac{1}{3}, 0, \frac{1}{5}, \ldots\right) \oplus \left(0, \frac{1}{2}, 0, \frac{1}{4}, 0, \ldots\right).
\]
Note that $$\sum_{k=1}^{\infty} \left|\frac{1}{2k-1}\right|^2 (2k-1)^2= \infty \text{~and~} \sum_{k=1}^{\infty} \left|\frac{1}{2k}\right|^2 (2k-1)^2= \infty.$$
This shows neither of these components of \( \Big(\frac{1}{n}\Big) \) belongs to $\mathcal{D}(T)$. Since \( \Big(\frac{1}{n}\Big) \) has a unique decomposition with respect to \( M \oplus M\p \), it follows that \( \mathcal{D}(T) \not\subseteq (M \cap \mathcal{D}(T)) \oplus (M\p \cap \mathcal{D}(T)) \).
Thus, \( M \) is not \( T \)-decomposable.
\end{example}

\noindent We will provide one more example of a subspace in the function space $L_2[0,1]$ that is not $T$-decomposable.
\begin{example}
    Suppose that $\mathcal{H} = L_2[0,1]$ and $M = \{f \in \mathcal{H} : f \equiv 0~ \text{on}~ [0,\frac{1}{2}]\}$.\\
    Let $T$ be an operator on $\mathcal{H}$ with $\mathcal{D}(T) = \{\text{polynomials on } [0,1]\}$.

\noindent Consider $N = \{f \in \mathcal{H} : f \equiv 0~ \text{on}~ \big[\frac{1}{2},1\big]\}$. Suppose $g \in N$. Then $fg \equiv 0$ for all $f \in M$. So, 
$$
\langle f, g \rangle = 0 \quad \text{for all } f \in M.
$$
Thus, 
$$
N \subseteq M\p.
$$
Suppose $g \notin N$. Then there exists $(a,b) \subseteq \big[\frac{1}{2},1\big]$ such that either $g > 0$ or $g < 0$ on $(a,b)$. Let $f \in \mathcal{H}$ such that $f = 0$ on $\big[0,a\big] \cup \big[b,1\big]$ and $f = 1$ on $(a,b)$. Now, either $fg > 0$ or $fg < 0$ on $(a,b)$ and $fg = 0$ on $\big[0,a\big] \cup \big[b,1\big]$. So, 
$$
\langle f, g \rangle \neq 0.
$$
Thus, 
$$
M\p \subseteq N.
$$
Therefore, 
$$
M\p = \Big\{f \in \mathcal{H} : f \equiv 0~ \text{on}~ \Big[\frac{1}{2},1\Big]\Big\}.
$$
Now, suppose $f \in \mathcal{D}(T) \cap M$. Then $f \in M$, which means $f$ has uncountably many zeros. Since $f \in \mathcal{D}(T)$, $f$ is a polynomial. However, a polynomial can only have finitely many zeros, leading to a contradiction. Hence, 
$$
\mathcal{D}(T) \cap M = \{0\}.
$$
Suppose $f \in \mathcal{D}(T) \cap M\p$. Then $f \in M\p$, which implies $f$ has uncountably many zeros. Again, since $f \in \mathcal{D}(T)$, $f$ is a polynomial, which can only have finitely many zeros, this is a contradiction. Hence, 
$$
\mathcal{D}(T) \cap M\p = \{0\}.
$$
\noindent Thus, 
$$
\mathcal{D}(T) \cap M \oplus \mathcal{D}(T) \cap M\p = \{0\} \oplus \{0\} = \{0\} \neq \mathcal{D}(T).
$$
Hence, $M$ is not $T$-decomposable.
\end{example}


A natural question arises: For a given operator \(T\) with \(\mathcal{D}(T) \neq \mathcal{H}\), does there exist a subspace that is not \(T\)-decomposable? The following theorem answers this question for the case when \(T\) is a densely defined operator.

\begin{theorem}
Let $T$ be a proper densely defined operator on $\mathcal{H}$. Then there exists a one-dimensional subspace $M$ of $\mathcal{H}$ such that $M$ is not $T$-decomposable.
\end{theorem}

\begin{proof}
Since $\mathcal{D}(T)$ is not closed, there exists $x_0 \in \overline{\mathcal{D}(T)} \setminus \mathcal{D}(T)$ with $\|x_0\| = 1$. Let $\mathcal{A}$ be an maximal orthonormal set in $\overline{\mathcal{D}(T)}$ containing $x_0$. Define 
$$
\mathcal{A}_1 = \mathcal{A} \cap \mathcal{D}(T), \quad \mathcal{A}_2 = \mathcal{A} \cap (\mathcal{H} \setminus \mathcal{D}(T)).
$$
Thus, by construction, we have $x_0 \in \mathcal{A}_2$ and $x_0 \in \left( \overline{\operatorname{span}(\mathcal{A}_1)} \right)^\perp $.

\noindent \textbf{Claim:} There exists $x \in \mathcal{D}(T)$ such that $x = y + z$, where $y \in \overline{\operatorname{span}(\mathcal{A}_1)}$, $z \in \overline{\operatorname{span}(\mathcal{A}_2)}$ and $z \neq 0$.

\noindent Assume, for contradiction, that the claim is false. Then we must have
$$
\mathcal{D}(T) \subseteq \overline{\operatorname{span}(\mathcal{A}_1)}.
$$
As a result,
$$
x_0 \in \left( \overline{\operatorname{span}(\mathcal{A}_1)} \right)^\perp 
\subseteq \left( \mathcal{D}(T) \right)^\perp 
= \left( \overline{\mathcal{D}(T)} \right)^\perp = \{0\}.
$$
This implies that $x_0 = 0$, which contradicts our assumption.
\noindent Hence, the original claim is true.

\noindent Now, write 
$$
y = \sum_{y_i \in \mathcal{A}_1} \alpha_i y_i, \quad z = \sum_{z_i \in \mathcal{A}_2} \beta_i z_i.
$$
Since $z \neq 0$, there exists some $i_0$ such that $\beta_{i_0} \neq 0$. Let $M = \operatorname{span}(z_{i_0})$. Then 
$$
M\p = (\operatorname{span}(z_{i_0}))\p = \operatorname{span}(\mathcal{A} \setminus \{z_{i_0}\}).
$$
Now, 
$$
x = y + z = \sum_{y_i \in \mathcal{A}_1} \alpha_i y_i + \sum_{z_i \in \mathcal{A}_2, z_i \neq z_{i_0}} \beta_i z_i + \beta_{i_0} z_{i_0}=\sum_{w_i \in \mathcal{A}_1 \cup \mathcal{A}_2, w_i \neq z_{i_0}} \gamma_i w_i +  \beta_{i_0} z_{i_0}.
$$
Define 
$$
a = \sum_{w_i \in \mathcal{A}_1 \cup \mathcal{A}_2, w_i \neq z_{i_0}} \gamma_i w_i, \quad b = \beta_{i_0} z_{i_0},
$$
where $a \in M\p$ and $b \in M$. Since $z_{i_0} \notin \mathcal{D}(T)$, we have $z_{i_0} \notin \mathcal{D}(T) \cap M$. Thus, $a \notin \mathcal{D}(T) \cap M$ and $b \notin \mathcal{D}(T) \cap M\p$.

Since this sum is unique with respect to $M \oplus M\p$, we cannot write $x = c + d,$ for some $c \in \mathcal{D}(T) \cap M$ and $d \in \mathcal{D}(T) \cap M\p$. Hence, 
$$
\mathcal{D}(T) \neq (\mathcal{D}(T) \cap M) \oplus (\mathcal{D}(T) \cap M\p).
$$
Therefore, $M$ is not $T$-decomposable.
\end{proof}
Now, we provide some characterizations for $T$-decomposable subspaces in terms of orthogonal projections onto the subspace $M$.
  \begin{theorem}\label{tdecchar}
        Let $T: \mathcal{D}(T) \subseteq \mathcal{H} \to \mathcal{K}$ and $M$ be a closed subspace of $\mathcal{H}$. Then $M$ is $T$-decomposable if and only if $P_{M}(\mathcal{D}(T)) \subseteq \mathcal{D}(T).$
    \end{theorem}

\begin{proof}
    Assume that $M$ is $T$-decomposable. Then we have $$\mathcal{D}(T) = (\mathcal{D}(T) \cap M) \oplus (\mathcal{D}(T) \cap M\p).$$
    Therefore any $x \in \mathcal{D}(T) $ can be written as $x=y+z,$ for some $y \in \mathcal{D}(T) \cap M$ and $z \in \mathcal{D}(T) \cap M\p.$ Hence we get that $P_M x= y \in \mathcal{D}(T) \cap M.$ Thus $$P_{M}(\mathcal{D}(T)) \subseteq \mathcal{D}(T).$$

   Conversely, assume that $P_M(\mathcal{D}(T)) \subseteq \mathcal{D}(T).$ Since $\mathcal{H} = M \oplus M^\perp$, any $x \in \mathcal{D}(T)$ can be uniquely expressed as $x = y + z$, where $y \in M$ and $z \in M^\perp.$ 
Applying $P_M$, we get $P_M x = y.$ By assumption, $P_M x \in \mathcal{D}(T)$, so $y \in \mathcal{D}(T).$ Since $z = x - y$ and both $x$ and $y$ belong to $\mathcal{D}(T)$, it follows that $z \in \mathcal{D}(T).$ 
Hence, $y \in \mathcal{D}(T) \cap M$ and $z \in \mathcal{D}(T) \cap M^\perp.$
    Thus $$\mathcal{D}(T) \subseteq (\mathcal{D}(T) \cap M) \oplus (\mathcal{D}(T) \cap M\p).$$
    Hence $M$ is $T$-decomposable.
\end{proof}
    \begin{Lemma}\label{tdeclem}
            Let $T: \mathcal{D}(T) \subseteq \mathcal{H} \to \mathcal{K}$ and $M$ be a closed subspace of $\mathcal{H}$. Then  $P_{M}(\mathcal{D}(T)) \subseteq \mathcal{D}(T)$ if and only if $P_{M\p}(\mathcal{D}(T)) \subseteq \mathcal{D}(T)$.
    \end{Lemma}
    \begin{proof}
        Assume that $P_{M}(\mathcal{D}(T)) \subseteq \mathcal{D}(T).$ Since $\mathcal{H}= M \oplus M\p,$ we can write any $x \in \mathcal{D}(T)$ as $x =y+z,$ for some $y \in M$ and $z \in M\p.$ Thus, we have $P_M x= y.$ Since we have assumed that $P_{M}(\mathcal{D}(T)) \subseteq \mathcal{D}(T)$, we have $P_M x= y \in \mathcal{D}(T).$ From the above facts, we also conclude that $z=x-y \in \mathcal{D}(T).$ Thus we have $z \in \mathcal{D}(T) \cap M\p.$
        Therefore $P_{M\p} x = z \in \mathcal{D}(T) \cap M\p.$ Hence $$P_{M\p} (\mathcal{D}(T)) \subseteq \mathcal{D}(T) .$$

       \noindent Converse part can be proved by replacing $M$ by $M\p$.
    \end{proof}
    \begin{corollary}
          Let $T: \mathcal{D}(T) \subseteq \mathcal{H} \to \mathcal{K}$ and $M$ be a closed subspace of $\mathcal{H}$. Then  $M$ is $T$-decomposable if and only if $M\p$ is $T$-decomposable.
    \end{corollary}
\begin{proof}
By Theorem \ref{tdecchar}, $M$ is $T$-decomposable if and only if 
$$
P_{M}(\mathcal{D}(T)) \subseteq \mathcal{D}(T).
$$
Now, by Lemma \ref{tdeclem}, we have $P_{M\p}(\mathcal{D}(T)) \subseteq \mathcal{D}(T)$ if and only if $P_{M}(\mathcal{D}(T)) \subseteq \mathcal{D}(T)$. Combining both results, we have $M$ is $T$-decomposable if and only if 
$$
P_{M\p}(\mathcal{D}(T)) \subseteq \mathcal{D}(T).
$$
Again, by Theorem \ref{tdecchar}, we have $M\p$ is $T$-decomposable if and only if 
$$
P_{M\p}(\mathcal{D}(T)) \subseteq \mathcal{D}(T).
$$
Combining the previous results, we conclude that $M$ is $T$-decomposable if and only if $M\p$ is $T$-decomposable.
\end{proof}

\noindent Let $T:\mathcal{D}(T) \subseteq \mathcal{H} \to \mathcal{K}$ with $\mathcal{D}(T)$ being a proper dense subspace of $\mathcal{H}$ and $M$ is a closed $T$-decomposable subspace of $\mathcal{H}$ and $N$ is a closed subspace of $\mathcal{K}$ respectively. We can express $\mathcal{D}(T)$ as the orthogonal direct sum of $ M_1 $ and $ M_2 $ and $\mathcal{K}$ as the orthogonal direct sum of $ N $ and $ N^\perp $. With respect to the above decompositions, the operator $ T $ can be written in a block matrix form as
\begin{equation}\label{*}
    T = \begin{pmatrix}
    A & B \\
    C & D
    \end{pmatrix}.
\end{equation}
From this notion, we have  
$
\mathcal{D}(A) = \mathcal{D}(C) = M_1  \text{ and }  \mathcal{D}(B) = \mathcal{D}(D) = M_2.
$\\
Therefore, if $T$ is written in the form (\ref{*}), then $$\mathcal{D}(T)=\mathcal{D}(A) \oplus \mathcal{D}(D)= \mathcal{D}(A) \oplus \mathcal{D}(B)= \mathcal{D}(C) \oplus \mathcal{D}(D)= \mathcal{D}(C) \oplus \mathcal{D}(B).$$

\noindent We know that for a bounded operator 
$$
T = \begin{pmatrix}
    A & B \\
    C & D
\end{pmatrix},
$$ 
the adjoint is given as 
$$
T^* = \begin{pmatrix}
    A^* & C^* \\
    B^* & D^*
\end{pmatrix}.
$$ 
However, this is not true in the case of a densely defined operator. The following theorem provides the criteria for when it holds.

\begin{theorem}
    Let $T:\mathcal{D}(T) \subseteq \mathcal{H} \to \mathcal{K}$ be a densely defined operator as defined in (\ref{*}). Define 
    $$
    T\q = \begin{pmatrix}
        A^* & C^* \\
        B^* & D^*
    \end{pmatrix}: \{\mathcal{D}(A^*) \cap \mathcal{D}(B^*)\} \oplus \{\mathcal{D}(C^*) \cap \mathcal{D}(D^*)\} \to \mathcal{H}.
    $$
    Then
 \begin{enumerate}
\item $T\q \subseteq T^*$.
\item $T\q = T^*$ if and only if $N$ is $T^*$-decomposable.
\end{enumerate}
\end{theorem}
\begin{proof}
The domain of $T^*$ is given by:
$$
\mathcal{D}(T^*) := \left\{ y \in \mathcal{K} : \exists\, z_y \in \mathcal{H} \text{ such that } \langle Tx, y \rangle = \langle x, z_y \rangle \ \forall x \in \mathcal{D}(T) \right\}.
$$

Now for any $x = x_1 \oplus x_2 \in \mathcal{D}(T)$ and $y = y_1 \oplus y_2 \in \mathcal{D}(T\q)$,
\begin{align*}
\langle Tx, y \rangle &= \langle (Ax_1 + Bx_2) \oplus (Cx_1 + Dx_2), y_1 \oplus y_2 \rangle \\
&= \langle Ax_1, y_1 \rangle + \langle Cx_1, y_2 \rangle + \langle Bx_2, y_1 \rangle + \langle Dx_2, y_2 \rangle \\
&= \langle x_1, A^* y_1 \rangle + \langle x_1, C^* y_2 \rangle + \langle x_2, B^* y_1 \rangle + \langle x_2, D^* y_2 \rangle \\
&= \langle x_1 \oplus x_2, A^* y_1 \oplus 0 \rangle + \langle x_1 \oplus x_2, C^* y_2 \oplus 0 \rangle \\
&~~~+ \langle x_1 \oplus x_2, 0 \oplus B^* y_1 \rangle + \langle x_1 \oplus x_2, 0 \oplus D^* y_2 \rangle \\
&= \langle x_1 \oplus x_2, (A^* + B^*) y_1 \oplus (C^* + D^*) y_2 \rangle \\
&= \langle x, T\q y \rangle.
\end{align*}
Therefore, $y \in \mathcal{D}(T^*)$.  
Thus $T\q \subseteq T\s.$

Conversely, assume that $N$ is $T^*$-decomposable and with respect to $N \oplus N\p,$ $T\s$ has the following form:
$$
T^* = \begin{pmatrix} A_1^* & C_1^* \\ B_1^* & D_1^* \end{pmatrix},
$$
with
$$
\mathcal{D}(T^*) =  \mathcal{D}(A_1^*)  \oplus  \mathcal{D}(D_1^*), ~~\mathcal{D}(A_1^*) =\mathcal{D}(B_1^*) \text{ and } \mathcal{D}(C_1^*)=\mathcal{D}(D_1^*).
$$

Now for any $x = x_1 \oplus x_2 \in \mathcal{D}(T)$ and $y = y_1 \oplus y_2 \in \mathcal{D}(T^*)$, where $x_1 \in M_1$, $x_2 \in M_2$, $y_1 \in \mathcal{D}(A_1^*)$, and $y_2 \in \mathcal{D}(D_1^*)$, we have
\begin{equation}
\langle Tx, y \rangle = \langle x, T^* y \rangle. 
\end{equation}
Take $x_2 = 0$ and $y_2 = 0$. Then for $x_1 \oplus 0 \in M_1$ and $y_1 \oplus 0 \in \mathcal{D}(T^*)$, we get
\begin{align*}
\langle T(x_1 \oplus 0), y_1 \oplus 0 \rangle &= \langle x_1 \oplus 0, T^*(y_1 \oplus 0) \rangle \\
\Rightarrow \langle Ax_1 \oplus Cx_1, y_1 \oplus 0 \rangle &= \langle x_1 \oplus 0, A_1^* y_1 \oplus B_1^* y_1 \rangle \\
\Rightarrow \langle Ax_1, y_1 \rangle &= \langle x_1, A_1^* y_1 \rangle.
\end{align*}
Hence, $y_1 \in \mathcal{D}(A^*)$.  
Similarly, by taking $x_1 = 0$ and $y_2 = 0$, we get $y_1 \in \mathcal{D}(B^*)$. Thus,
$$
\mathcal{D}(A_1^*) = \mathcal{D}(B_1^*) \subseteq \mathcal{D}(A^*) \cap \mathcal{D}(B^*).
$$
Similarly, by taking $x_1 = 0$, $y_1 = 0$ and $x_2 = 0$, $y_1 = 0$, we get
$$
\mathcal{D}(C_1^*) = \mathcal{D}(D_1^*) \subseteq \mathcal{D}(C^*) \cap \mathcal{D}(D^*).
$$
Thus,
$$
\mathcal{D}(T^*) \subseteq \mathcal{D}(T\q).
$$
Since (1) holds, we conclude that $T^* = T\q$.
\end{proof}

\noindent The operator $T\q = \begin{pmatrix}
        A^* & C^* \\
        B^* & D^*
    \end{pmatrix}$ is called formal adjoint of $T.$\\
 Now onwards, we denote 
\[
N_1 := \mathcal{D}(A\s) \cap \mathcal{D}(B\s) 
\quad \text{and} \quad 
N_2 := \mathcal{D}(C\s) \cap \mathcal{D}(D\s).
\]
Throughout this paper, we consider $A\s$ and $B\s$ as the restrictions of the operators 
$A\s$ and $B\s$ to $N_1$, respectively, and $C\s$ and $D\s$ as the restrictions of 
$C\s$ and $D\s$ to $N_2$. With this notation, we have
\[
T\q = 
\begin{pmatrix}
A\s & C\s \\
B\s & D\s
\end{pmatrix}
\text{ with }
\mathcal{D}(T\q) = N_1 \oplus N_2,
\]
with $\mathcal{D}(A\s) = \mathcal{D}(B\s) = N_1$ and 
$\mathcal{D}(C\s) = \mathcal{D}(D\s) = N_2.$
 
\section{Densely defined complementable operators}

The concept of bounded complementable operators, as introduced by \cite{Antezana}, provides a framework for analyzing operators in terms of subspace decompositions and range inclusions. This idea is particularly powerful when operators can be expressed in block matrix form, allowing the study of individual components with respect to specific subspaces.

For unbounded operators, such decompositions and range considerations are not always available due to domain restrictions. Nevertheless, when an operator can be written in a block operator matrix form, we extend the characterizations of complementability given in \cite{Arias} for the triplets $(T,M,N)$ to densely defined possibly unbounded operators.
 \begin{definition}
    Let $T:\mathcal{D}(T) \subseteq \mathcal{H} \to \mathcal{K}$ be a densely defined operator as defined in (\ref{*}). Then $T$ is called $(M,N)$-complementable if $\mathcal{R}(C) \subseteq \mathcal{R}(D)$ and $\mathcal{R}(B\s) \subseteq \mathcal{R}(D\s)$.
 \end{definition}
 
If $T$ is of the form \eqref{*} and is $(M, N)$-complementable, then $\mathcal{R}(C) \subseteq \mathcal{R}(D)$ and $\mathcal{R}(B^*) \subseteq \mathcal{R}(D^*)$. Therefore, by Theorem \ref{Douglasunbdd}, there exist operators
$$
X : \mathcal{D}(C) \to \mathcal{D}(D) \quad \text{and} \quad Y : \mathcal{D}(B^*)  \to \mathcal{D}(D^*)
$$
such that
$$
C = D X  \text{ and }  B^* = D^* Y.
$$
Throughout our discussion, we will consider the operators $X$ and $Y$ as defined above.

The following results establish fundamental characterizations and properties of unbounded $(M, N)$-complementable operators. Building upon the foundational frameworks developed in \cite{Antezana} and \cite{Arias}, these results extend the theory to accommodate a wider class of structural configurations and functional behaviors inherent to such operators.
 \begin{theorem}
Let $T:\mathcal{D}(T) \subseteq \mathcal{H} \to \mathcal{K}$ be a densely defined operator as defined in (\ref{*}). Then $T$ is $(M, N)$-complementable if and only if for any pair of projections $(P_r, P_l)$ with $\mathcal{R}(P_r) = M\p$ and $\mathcal{N}(P_l) = N$, there exist operators $M_r$ and $M_l$ such that the following hold:
\begin{enumerate}
    \item $P_r M_r = M_r$, $P_l T M_r = P_l T$ on $\mathcal{D}(T)$,
    \item $P_l^* M_l = M_l$, $P_r^* T\q M_l = P_r^* T\q$ on $\mathcal{D}(T\q)$.
\end{enumerate}
\end{theorem}
\begin{proof}
Let $(P_r, P_l)$ be a pair of projections such that $\mathcal{R}(P_r) = M^\perp$ and $\mathcal{N}(P_l) = N$. Since $\mathcal{R}(C) \subseteq \mathcal{R}(D)$ and $\mathcal{R}(B^*) \subseteq \mathcal{R}(D^*)$, there exist operators $X$ and $Y$ with $C = DX$ and $B^* = D^*Y$. We may represent $P_r$ and $P_l$ in block form as
\[
P_r = \begin{pmatrix} 0 & 0 \\ E & I \end{pmatrix}, 
\qquad 
P_l = \begin{pmatrix} 0 & F \\ 0 & I \end{pmatrix},
\]
and define
\[
M_r = \begin{pmatrix} 0 & 0 \\ X & I \end{pmatrix}, 
\qquad 
M_l = \begin{pmatrix} 0 & 0 \\ Y & I \end{pmatrix}.
\]

Since $\mathcal{D}(P_r) = \mathcal{H}$, we have $\mathcal{D}(P_rM_r) = \mathcal{D}(M_r)$, and direct computation shows
\[
P_r M_r 
= 
\begin{pmatrix} 0 & 0 \\ E & I \end{pmatrix}
\begin{pmatrix} 0 & 0 \\ X & I \end{pmatrix}
= 
\begin{pmatrix} 0 & 0 \\ X & I \end{pmatrix} 
= M_r.
\]
Thus $P_rM_r = M_r$.

Next, observe that $\mathcal{D}(P_l) = \mathcal{K}$, hence $\mathcal{D}(P_lT) = \mathcal{D}(T)$ and $\mathcal{D}(P_lTM_r) = \mathcal{D}(TM_r)$. Therefore,
\[
\mathcal{D}(P_lT) \supseteq \mathcal{D}(P_lTM_r) \cap \mathcal{D}(T).
\]
Let $u \in \mathcal{D}(T)$ with $u = x+y$ where $x \in M_1$ and $y \in M_2$. Then
\[
M_r u = \begin{pmatrix} 0 \\ Xx + y \end{pmatrix}.
\]
Since $\mathcal{R}(X) \subseteq \mathcal{D}(D) = \mathcal{D}(B)$, it follows that $M_r u \in \mathcal{D}(T)$, hence $TM_r u \in \mathcal{D}(P_l) = \mathcal{K}$. Consequently,
\[
\mathcal{D}(P_lT) \subseteq \mathcal{D}(P_lTM_r) \cap \mathcal{D}(T),
\]
and thus
\[
\mathcal{D}(P_lT) \cap \mathcal{D}(T) = \mathcal{D}(P_lTM_r) \cap \mathcal{D}(T).
\]
Since $\mathcal{D}(X) = \mathcal{D}(C)$, we also have $\mathcal{D}(T) \subseteq \mathcal{D}(M_r)$.

Now, for $u = x+y \in \mathcal{D}(T)$, a direct calculation gives
\begin{align*}
P_lTM_r(x+y) 
&= \begin{pmatrix} 0 & F \\ 0 & I \end{pmatrix}
\begin{pmatrix} A & B \\ C & D \end{pmatrix}
\begin{pmatrix} 0 & 0 \\ X & I \end{pmatrix}
\begin{pmatrix} x \\ y \end{pmatrix} \\
&= \begin{pmatrix} 0 & F \\ 0 & I \end{pmatrix}
\begin{pmatrix} B(Xx+y) \\ D(Xx+y) \end{pmatrix} \\
&= \begin{pmatrix} FCx + FDy \\ Cx + Dy \end{pmatrix} \\
&= \begin{pmatrix} 0 & F \\ 0 & I \end{pmatrix}
\begin{pmatrix} A & B \\ C & D \end{pmatrix}
\begin{pmatrix} x \\ y \end{pmatrix} \\
&= P_lTu.
\end{align*}
Thus $P_lTM_r = P_lT$ on $\mathcal{D}(T)$.

Similarly, on $\mathcal{D}(T\q)$, one checks that $P_r^\ast T\q M_l = P_r^\ast T\q$ and $P_l^\ast M_l = M_l$. This completes the proof of necessity.

\medskip
Conversely, let $(P_r,P_l)$ be a pair of projections with $\mathcal{R}(P_r) = M^\perp$ and $\mathcal{N}(P_l) = N$, and let $M_r, M_l$ be operators satisfying the given conditions. Writing $P_r$ and $P_l$ as before, and using $P_rM_r = M_r$ and $P_l^\ast M_l = M_l$, we may represent
\[
M_r = \begin{pmatrix} 0 & 0 \\ J_1 & K_1 \end{pmatrix}, 
\qquad 
M_l = \begin{pmatrix} 0 & 0 \\ J_2 & K_2 \end{pmatrix}.
\]

The conditions $P_lTM_r = P_lT$ and $P_r^\ast T^\ast M_l = P_r^\ast T^\ast$ then yield
$\mathcal{R}(C) \subseteq \mathcal{R}(D)$ and $\mathcal{R}(B^*) \subseteq \mathcal{R}(D^*)$ respectively.\\
Thus, $T$ is $(M, N)$-complementable.
\end{proof}

The complementability of an operator can also be characterized using the existence of projections that link its null spaces to subspaces. This perspective leads to the following equivalent condition for $(M,N)$-complementability, involving projections with prescribed null spaces and corresponding range inclusions.
\begin{theorem}
Let $T:\mathcal{D}(T) \subseteq \mathcal{H} \to \mathcal{K}$ be a densely defined operator as defined in (\ref{*}). Then, the operator $T$ is $(M, N)$-complementable if and only if there exist projections $P$ and $Q$ with $\mathcal{N}(P)=M\p$, $\mathcal{N}(Q)=N\p$ such that $\mathcal{R}(TP) \subseteq N$ on $\mathcal{D}(T)$ and $\mathcal{R}(T\q Q) \subseteq M$ on $\mathcal{D}(T\q)$.
\end{theorem}
\begin{proof}
Suppose first that $T$ is $(M,N)$-complementable. 
Consider the projections
\[
P = \begin{pmatrix} I & 0 \\ -X & 0 \end{pmatrix},
\qquad 
Q = \begin{pmatrix} I & 0 \\ -Y & 0 \end{pmatrix}.
\]
Let $u \in \mathcal{D}(T)$, and write $u = x+y$ with $x \in M_1$ and $y \in M_2$. Then
\[
TPu 
= 
\begin{pmatrix} A & B \\ C & D \end{pmatrix}
\begin{pmatrix} I & 0 \\ -X & 0 \end{pmatrix}
\begin{pmatrix} x \\ y \end{pmatrix}
=
\begin{pmatrix} A & B \\ C & D \end{pmatrix}
\begin{pmatrix} x \\ -Xx \end{pmatrix}
=
\begin{pmatrix} Ax - BXx \\ Cx - DXx \end{pmatrix}.
\]
Since $C = DX$, it follows that
\[
TPu 
= 
\begin{pmatrix} Ax - BXx \\ 0 \end{pmatrix}
=
\begin{pmatrix} A - BX & 0 \\ 0 & 0 \end{pmatrix}
\begin{pmatrix} x \\ y \end{pmatrix}.
\]
Hence
\[
TP = \begin{pmatrix} A - BX & 0 \\ 0 & 0 \end{pmatrix}
\quad \text{on } \mathcal{D}(T).
\]
Similarly, we compute
\[
T\q Q 
= 
\begin{pmatrix} A^* & C^* \\ B^* & D^* \end{pmatrix}
\begin{pmatrix} I & 0 \\ -Y & 0 \end{pmatrix}
= 
\begin{pmatrix} A^* - B^*Y & 0 \\ 0 & 0 \end{pmatrix}
\quad \text{on } \mathcal{D}(T\q).
\]
Since $\mathcal{R}(TP) \subseteq N$ and $\mathcal{R}(T\q Q) \subseteq M$, the required conditions are satisfied.

\medskip
Conversely, suppose there exist projections $P$ and $Q$ with $\mathcal{N}(P) = M^\perp$, $\mathcal{N}(Q) = N^\perp$ such that $\mathcal{R}(TP) \subseteq N$ and $\mathcal{R}(T^*Q) \subseteq M$. We may write
\[
P = \begin{pmatrix} I & 0 \\ E & 0 \end{pmatrix},
\qquad
Q = \begin{pmatrix} I & 0 \\ F & 0 \end{pmatrix},
\]
where $\mathcal{R}(E) \subseteq M_2$ and $\mathcal{R}(F) \subseteq N_2$.
Then
\[
TP 
= 
\begin{pmatrix} A & B \\ C & D \end{pmatrix}
\begin{pmatrix} I & 0 \\ E & 0 \end{pmatrix}
=
\begin{pmatrix} A + BE & 0 \\ C + DE & 0 \end{pmatrix},
\quad \text{on } \mathcal{D}(T),
\]
and
\[
T^*Q
=
\begin{pmatrix} A^* & C^* \\ B^* & D^* \end{pmatrix}
\begin{pmatrix} I & 0 \\ F & 0 \end{pmatrix}
=
\begin{pmatrix} A^* + C^*F & 0 \\ B^* + D^*F & 0 \end{pmatrix},
\quad \text{on } \mathcal{D}(T^*).
\]
Since $\mathcal{R}(TP) \subseteq N$ and $\mathcal{R}(T^*Q) \subseteq M$, it follows that
$
C = -DE
\text{ and } 
B^* = -D^*F.$
Consequently,
\[
\mathcal{R}(C) \subseteq \mathcal{R}(D),
\qquad
\mathcal{R}(B^*) \subseteq \mathcal{R}(D^*).
\]

Therefore, $T$ is $(M,N)$-complementable.
\end{proof}

Beyond range and projection-based characterizations, the concept of $(M, N)$-complementability can also be understood through the lens of domain decomposition. This perspective emphasizes how the domain of the operator and its adjoint relate to the underlying subspaces $M$ and $N$.
By expressing the domains of $T$ and $T\q$ as sums of inverse images and complementary components, one obtains a structural insight into how the operator interacts with the subspaces. The following theorems express $(M, N)$-complementability in terms of domain decomposition, providing a direct relationship between the operator's domain and the involved subspaces:
\begin{theorem}\label{t319}
Let $T:\mathcal{D}(T) \subseteq \mathcal{H} \to \mathcal{K}$ be a densely defined operator as defined in (\ref{*}). Then $T$ is $(M,N)$-complementable if and only if 
$$\mathcal{D}(T) = T^{-1}(N) + M_2 \quad \text{and} \quad \mathcal{D}(T\q) = (T\q)^{-1}(M) + N_2.$$
\end{theorem}

\begin{proof} 
Let $x \in \mathcal{D}(T)$. We can decompose $x$ as $x = y + z$, where $y \in M_1$ and $z \in M_2$. Since $-Xy \in \mathcal{D}(D) = M_2$, we have $y - Xy \in \mathcal{D}(T)$.\\
Next, we have
$$
T(y - Xy) = \begin{pmatrix} A & B \\ C & D \end{pmatrix} \begin{pmatrix} y \\ -Xy \end{pmatrix} = \begin{pmatrix} (A - BX)y \\ (C - DX)y \end{pmatrix} = \begin{pmatrix} (A - BX)y \\ 0 \end{pmatrix} \in N.
$$
Thus, $y - Xy \in T^{-1}(N)$. Therefore, we can express $$x = y + z = y - Xy + z + Xy \in T^{-1}(N) + M_2.$$ Clearly, $T^{-1}(N), M_2 \subseteq \mathcal{D}(T)$ and hence we conclude that $$T^{-1}(N) + M_2 = \mathcal{D}(T).$$

Now, let $x \in \mathcal{D}(T\q)$. We write $x = y + z$, where $y \in N_1$ and $z \in N_2$. Since $-Yy \in \mathcal{D}(D\s) = N_2$, it follows that $y - Yy \in \mathcal{D}(T\q)$.\\
Again, we have
$$
T\q(y - Yy) = \begin{pmatrix} A\s & C\s \\ B\s & D\s \end{pmatrix} \begin{pmatrix} y \\ -Yy \end{pmatrix} = \begin{pmatrix} (A\s - C\s Y)y \\ (B\s - D\s Y)y \end{pmatrix} = \begin{pmatrix} (A\s - C\s Y)y \\ 0 \end{pmatrix} \in M.
$$
Thus, $y - Yy \in (T\q)^{-1}(M)$. Therefore, we can express $$x = y + z = y - Yy + z + Yy \in (T\q)^{-1}(M) + N_2.$$ Again, $(T\q)^{-1}(M), N_2 \subseteq \mathcal{D}(T\q)$, so we conclude that $$(T\q)^{-1}(M) + N_2 = \mathcal{D}(T\q).$$
  
Conversely, suppose $x \in M_1$. By assumption, there exist elements $y$ and $z$ such that $y \in M_2$ and $Tz \in N$ with $x = y + z$. Now,
$$
Tz = \begin{pmatrix} A & B \\ C & D \end{pmatrix} \begin{pmatrix} x \\ -y \end{pmatrix} = \begin{pmatrix} Ax - By \\ Cx - Dy \end{pmatrix}.
$$
Since $Tz \in N$, it follows that $Cx = Dy$ and hence $\mathcal{R}(C) \subseteq \mathcal{R}(D)$.

Next, let $x \in N_1$. There exist $y$ and $z$ such that $y \in N_2$ and $T\q z \in M$ with $x = y + z$. We have
$$
T\q z = \begin{pmatrix} A\s & C\s \\ B\s & D\s \end{pmatrix} \begin{pmatrix} x \\ -y \end{pmatrix} = \begin{pmatrix} A\s x - C\s y \\ B\s x - D\s y \end{pmatrix}.
$$
Since $T\q z \in M$, it follows that $B\s x = D\s y$ and hence $\mathcal{R}(B\s) \subseteq \mathcal{R}(D\s)$.

\noindent Therefore, $T$ is $(M, N)$-complementable. 
\end{proof}

The following result provides an equivalent formulation of complementability in terms of domain decompositions involving dense subsets of $N\p$ and $M\p$. It refines earlier decompositions by expressing them through intersections with the domain and orthogonal complements of image subspaces, thus offering a more operational perspective that is especially useful in functional analysis and operator theory.

  \begin{theorem}
    Let $T:\mathcal{D}(T) \subseteq \mathcal{H} \to \mathcal{K}$ be a densely defined operator as defined in (\ref{*}). Let $N_2$ be dense in $N\p$. Then $T$ is $(M,N)$-complementable if and only if $$\mathcal{D}(T){=}M_2+((T\q(N_2))\p \cap \mathcal{D}(T)) \text{~and~} \mathcal{D}(T\q){=}N_2+((T(M_2))\p \cap \mathcal{D}(T\q)).$$
 \end{theorem}

\begin{proof}
    By Theorem \ref{t319}, it is sufficient to show that
    $$ T^{-1}(N) = (T\q(N_2))\p \cap \mathcal{D}(T)  \text{~and~}  (T\q)^{-1}(M) = (T(M_2))\p \cap \mathcal{D}(T\q). $$
  Let $x \in T^{-1}(N)$. Then, $$x \in \mathcal{D}(T) \text{~and~} \langle T x, m \rangle = 0 \text{~for all~} m \in N\p.$$
  This implies that $$\langle x, T\s m \rangle = 0 \text{~for all~} m \in N\p \cap \mathcal{D}(T\s).$$
  Therefore, $$x \in (N\p \cap \mathcal{D}(T\s))\p \subseteq (T\q(N_2))\p.$$
Hence, we conclude that
    $$ T^{-1}(N) \subseteq (T\q(N_2))\p \cap \mathcal{D}(T). $$
Now, we prove the reverse inclusion. Let $x \in (T\q(N_2))\p \cap \mathcal{D}(T)$. Since $\langle x, T\q m \rangle = 0$ for all $m \in N_2$ and $T\q \subseteq T\s$, we have $\langle x, T\s m \rangle = 0$ for all $m \in N_2$. Consequently, $$\langle T x, m \rangle = 0 \text{ for all } m \in N_2.$$ Since $N_2$ is dense in $N\p$, it follows that $$\langle T x, m \rangle = 0 \text{~for all~} m \in N\p.$$ Thus, $x \in T^{-1}(N)$, which gives
    $$ T^{-1}(N) \supseteq (T\q(N_2))\p \cap \mathcal{D}(T). $$
 Therefore, we have
    $$ T^{-1}(N) = (T\q(N_2))\p \cap \mathcal{D}(T). $$
    
  To prove the other equality, let us consider $x \in (T\q)^{-1}(M)$.
  Then, clearly, $x \in \mathcal{D}(T\q)$. Since $\langle T\q x, m \rangle = 0$ for all $m \in M\p$ and $T\q \subseteq T\s$, we also have $$\langle x, {{T\s} \s} m \rangle = 0 \text{ for all } m \in M\p \cap \mathcal{D}({{T\s} \s}).$$ Since ${T\s} \s = T$ on $\mathcal{D}(T)$, it follows that $$\langle x, T m \rangle = 0 \text{ for all } m \in M\p \cap \mathcal{D}(T).$$ Thus, $x \in (T(M_2))\p$. Therefore, we conclude that
    $$ (T\q)^{-1}(M) \subseteq (T(M_2))\p \cap \mathcal{D}(T\q). $$
   Now, we prove the reverse inclusion. Let $x \in (T(M_2))\p \cap \mathcal{D}(T\q)$. Since $\langle x, T m \rangle = 0$ for all $m \in M_2$, we also have $$\langle T\q x, m \rangle = 0 \text{ for all } m \in M_2.$$ Since $M_2$ is dense in $M\p$, it follows that $$\langle T\q x, m \rangle = 0 \text{~for all~} m \in M\p.$$ Thus, $x \in (T\q)^{-1}(M)$, which gives
    $$ (T\q)^{-1}(M) \supseteq (T(M_2))\p \cap \mathcal{D}(T\q). $$
 Therefore, we conclude that
    $ (T\q)^{-1}(M) = (T(M_2))\p \cap \mathcal{D}(T\q).$
\end{proof}

	Complementability can also be characterized in terms of the range of the operator and its adjoint. This perspective reveals how the image of the operator aligns with the given subspaces and complements the earlier domain-based formulations.

The decomposition of the range into components contributed by the subspaces $M$ and $N$ and their intersections with the range provides a geometric understanding of how the operator projects into the target subspaces. The following result captures this viewpoint.
 \begin{theorem}
Let $T:\mathcal{D}(T) \subseteq \mathcal{H} \to \mathcal{K}$ be a densely defined operator as defined in (\ref{*}). Then $T$ is $(M,N)$-complementable if and only if $\mathcal{R}(T)=T(M_2)+(N \cap \mathcal{R}(T))$ and $\mathcal{R}(T\q)=T\q(N_2)+(M\cap \mathcal{R}(T\q))$.
 \end{theorem}

\begin{proof}
    From Theorem \ref{t319}, we know that $T$ is $(M,N)$-complementable if and only if 
    $$
    \mathcal{D}(T) = T^{-1}(N) + M_2 \quad \text{and} \quad \mathcal{D}(T\q) = (T\q)^{-1}(M) + N_2.
    $$
    Thus, it suffices to prove the following equivalences:
    \begin{itemize}
        \item[(i)] $\mathcal{D}(T) = T^{-1}(N) + M_2$ if and only if $\mathcal{R}(T) = T(M_2) + (N \cap \mathcal{R}(T))$.
        \item[(ii)] $\mathcal{D}(T\q) = (T\q)^{-1}(M) + N_2$ if and only if $\mathcal{R}(T\q) = T\q(N_2) + (M \cap \mathcal{R}(T\q))$.
    \end{itemize}
    We first prove part (i): $\mathcal{D}(T) = T^{-1}(N) + M_2$ if and only if $\mathcal{R}(T) = T(M_2) + (N \cap \mathcal{R}(T))$.\\
$\implies :$ Assume $\mathcal{D}(T) = T^{-1}(N) + M_2$. Then,
    $$
    \mathcal{R}(T) = T(\mathcal{D}(T)) = T(T^{-1}(N) + M_2)= T(T^{-1}(N)) + T(M_2).
    $$
    Since $T(T^{-1}(N)) = N \cap \mathcal{R}(T)$, we have
    $$
    \mathcal{R}(T) = (N \cap \mathcal{R}(T)) + T(M_2).
    $$
    
 \noindent $\impliedby :$ Assume $\mathcal{R}(T) = T(M_2) + (N \cap \mathcal{R}(T))$. Then,
    $$
    \mathcal{D}(T) = T^{-1}(\mathcal{R}(T)) = T^{-1}(T(M_2) + (N \cap \mathcal{R}(T))).
    $$
    Using the distributive property of preimages under $T$, we get
    $$
    \mathcal{D}(T) = T^{-1}(T(M_2)) + T^{-1}(N \cap \mathcal{R}(T)).
    $$
    Since $T^{-1}(T(M_2)) = M_2 + \mathcal{N}(T)$ and $T^{-1}(N \cap \mathcal{R}(T)) = T^{-1}(N)$, it follows that
    $$
    \mathcal{D}(T) = T^{-1}(N) + M_2 + \mathcal{N}(T).
    $$
    As $\mathcal{N}(T) \subseteq T^{-1}(N)$, we conclude that
    $$
    \mathcal{D}(T) = T^{-1}(N) + M_2.
    $$

    \noindent $(ii)$ The \, proof \, of\, the \, equivalence \, $\mathcal{D}(T\q) = (T\q)^{-1}(M) + N_2$ \,if \,and \,only \,if $\mathcal{R}(T\q) = T\q(N_2) + (M \cap \mathcal{R}(T\q))$ follows similarly by replacing $T$ with $T\q$ and $M_2$ with $N_2$.
\end{proof}

\begin{Lemma}\label{Lemmaunbdd}
Let $T:\mathcal{D}(T) \subseteq \mathcal{H} \to \mathcal{K}$ be a densely defined operator as defined in (\ref{*}). Let $T$ be $(M,N)$-complementable and $Y$ be such that $B\s=D\s Y$. Then $\mathcal{R}(D) \subseteq \mathcal{D}(Y\s)$ and $B=Y\s D$ on $\mathcal{D}(B)=\mathcal{D}(D).$
 \end{Lemma}
\begin{proof}
Suppose that $B^* = D^* Y$. Let $x \in \mathcal{D}(D)$. Then for all $y \in \mathcal{D}(B^*) = \mathcal{D}(Y)$, we have
$$
\langle B^* y, x \rangle = \langle D^* Y y, x \rangle = \langle Y y, D x \rangle.
$$
Hence,
$$
\langle y, B x \rangle = \langle Y y, D x \rangle,
$$
for all $x \in \mathcal{D}(D)$ and $y \in \mathcal{D}(Y)$. This identity implies that for each $x \in \mathcal{D}(D)$, the element $D x$ defines a bounded conjugate-linear functional on $\mathcal{D}(Y)$ with respect to the pairing $\langle Y y, \cdot \rangle$. Thus, $D x \in \mathcal{D}(Y^*)$, and therefore,
$$
\mathcal{R}(D) \subseteq \mathcal{D}(Y^*).
$$
Since $B^* = D^* Y$, it follows that
$$
B^{**} \supseteq Y^* D^{**}.
$$
Moreover, as $\mathcal{D}(B) = \mathcal{D}(D) \subseteq \mathcal{D}(B^{**}) \cap \mathcal{D}(D^{**})$, we conclude that
$$
B = Y^* D \quad \text{on } \mathcal{D}(B) = \mathcal{D}(D).
$$
\end{proof}

	Another insightful perspective on $(M, N)$-complementability emerges from examining how the operator interacts with projections onto the subspaces $M$ and $N$. In this context, a remarkable symmetry condition involving the operator and its adjoint provides a powerful equivalence.
By requiring that the action of $T$ after projecting onto $M$ coincides with the adjoint of the action of $T\q$ after projecting onto $N$, one obtains a concise and elegant characterization of complementability. This symmetry is formalized in the following result.
 \begin{theorem}
 Let $T:\mathcal{D}(T) \subseteq \mathcal{H} \to \mathcal{K}$ be a densely defined operator as defined in (\ref{*}). Then $T$ is $(M,N)$-complementable if and only if there exist two projections $P$ and $Q$ with $\mathcal{N}(P)=M\p$, $\mathcal{N}(Q)=N\p$ such that $TP=(T\q Q)\s$ on $\mathcal{D}(T)$.
 \end{theorem}
\begin{proof}
Let 
$$
P = \begin{pmatrix}
    I & 0 \\
    -X & 0
\end{pmatrix} \quad \text{and} \quad 
Q = \begin{pmatrix}
    I & 0 \\
    -Y & 0
\end{pmatrix}.
$$ 
Clearly, $\mathcal{N}(P) = M^\perp$ and $\mathcal{N}(Q) = N^\perp$. Now, for $u=x \oplus y \in \mathcal{D}(T), Pu=\begin{pmatrix}
    I & 0 \\
    -X & 0
\end{pmatrix}\begin{pmatrix}
    x \\
    y
\end{pmatrix}=\begin{pmatrix}
    x \\
    -Xx
\end{pmatrix} \in \mathcal{D}(T),$ as $-Xx \in \mathcal{D}(D)=\mathcal{D}(B).$ Thus,
$$
TP = \begin{pmatrix}
    A & B \\
    C & D
\end{pmatrix}
\begin{pmatrix}
    I & 0 \\
    -X & 0
\end{pmatrix}
= \begin{pmatrix}
    A - BX & 0 \\
    0 & 0
\end{pmatrix} \text{ on } \mathcal{D}(T),
$$ and for $v=x \oplus y \in \mathcal{D}(T\q), Qv=\begin{pmatrix}
    I & 0 \\
    -Y & 0
\end{pmatrix}\begin{pmatrix}
    x \\
    y
\end{pmatrix}=\begin{pmatrix}
    x \\
    -Yx
\end{pmatrix} \in \mathcal{D}(T\q),$ as $-Yx \in \mathcal{D}(D\s)=\mathcal{D}(C\s).$ Thus,
$$
T^\times Q = \begin{pmatrix}
    A^* & C^* \\
    B^* & D^*
\end{pmatrix}
\begin{pmatrix}
    I & 0 \\
    -Y & 0
\end{pmatrix}
= \begin{pmatrix}
    A^* - C^*Y & 0 \\
    0 & 0
\end{pmatrix} \text{ on } \mathcal{D}(T^\times).
$$
Therefore,
$$
(T^\times Q)^* \supseteq \begin{pmatrix}
    A^{**} - Y^* C^{**} & 0 \\
    0 & 0
\end{pmatrix}.
$$
Given that $\mathcal{R}(C) \subseteq \mathcal{R}(D)$ and, by Lemma \ref{Lemmaunbdd}, $\mathcal{R}(C) \subseteq \mathcal{D}(Y)$, we deduce that 
$$
M_1 \subseteq \mathcal{D}(A^{**}) \cap \mathcal{D}(Y^* C^{**}).
$$
Hence,
$$
A^{**} - Y^* C^{**} = A - Y^* C \quad \text{on } M_1.
$$
Moreover, since $Y^* C = Y^* D X = BX$ on $M_1$, it follows that
$$
A - Y^* C = A - BX \quad \text{on } M_1.
$$
This establishes that $TP = (T^\times Q)^*$ on $\mathcal{D}(T)$.

Conversely, let
$$
P = \begin{pmatrix}
    I & 0 \\
    E & 0
\end{pmatrix}, \quad 
Q = \begin{pmatrix}
    I & 0 \\
    F & 0
\end{pmatrix}.
$$
Then,
$$
TP = \begin{pmatrix}
    A & B \\
    C & D
\end{pmatrix}
\begin{pmatrix}
    I & 0 \\
    E & 0
\end{pmatrix}
= \begin{pmatrix}
    BE + A & 0 \\
    C + DE & 0
\end{pmatrix},
$$
and
$$
T^\times Q = \begin{pmatrix}
    A^* & C^* \\
    B^* & D^*
\end{pmatrix}
\begin{pmatrix}
    I & 0 \\
    F & 0
\end{pmatrix}
= \begin{pmatrix}
    A^* + C^*F & 0 \\
    B^* + D^*F & 0
\end{pmatrix}.
$$
Consequently,
$$
(T^\times Q)^* \supseteq \begin{pmatrix}
    A^{**} + F^* C^{**} & B^{**} + F^* D^{**} \\
    0 & 0
\end{pmatrix}.
$$
Since $TP = (T^\times Q)^*$ on $\mathcal{D}(T)$, we obtain
$$
C + DE = 0 \quad \text{on } M_1, \quad \text{and} \quad B^{**} + F^* D^{**} = 0 \quad \text{on } M_2.
$$
From the identity $C + DE = 0$, it follows that $\mathcal{R}(C) \subseteq \mathcal{R}(D)$.  
In addition, the relation $B^{**} = -F^* D^{**}$ on $M_2$ implies $B^* = -D^* F^{**}$. Since $F \subseteq F^{**}$, we conclude that $B^* = -D^* F$ on $\mathcal{D}(B^*)$.
Thus, 
$$
\mathcal{R}(B^*) \subseteq \mathcal{R}(D^*).
$$
\end{proof}

    \section{Acknowledgement}
    ~\newline
     The first author thanks the Manipal Institute of Technology, Manipal Academy of Higher Education, Manipal for the financial support. The present work of the second author was partially supported by
   Anusandhan National Research Foundation (ANRF), Department of Science and Technology, Government
of India (Reference Number: MTR/2023/000471) under the scheme “Mathematical
Research Impact Centric Support (MATRICS)”.

\end{document}